\magnification =\magstep1

\baselineskip =13pt
\overfullrule =0pt

\vskip 3cm

 \centerline { \bf  THE SPACE OF 3-MANIFOLDS AND VASSILIEV}
\centerline {\bf  FINITE-TYPE INVARIANTS}

 \vskip .5cm

 \centerline { \bf  Nadya Shirokova\footnote{*}{\rm Department of Mathematics,
University of Chicago, Chicago IL 60637, \hfill\break email:
nadya@math.uchicago.edu}}

\vskip 1cm


\noindent { \bf  Abstract.}  In this paper we develop the theory of
finite-type invariants for homologically nontrivial 3-manifolds .
 We construct an infinite-dimensional affine space with
a hypersurface in it corresponding to manifolds with Morse singularities.
Connected components of the complement of this hypersurface correspond to
homeomorphism type of spin 3-manifolds. This suggests the natural axiomatics
of Vassiliev finite-type invariants for arbitrary closed 3-manifolds. An
example of  an invariant of order 1 is given.

\vskip 1cm

\centerline  { \bf  0. Introduction.}

\vskip 1cm

The idea which goes back to V. Arnold is to complete the space of all
objects by the degenerate ones so as to get a Euclidean space and then
study the topology of the degenerate locus which is related to the topology
of its complement via Alexander duality. This approach had a beautiful
application in the theory of knots. In 1986 V.Vassiliev constructed
so-called finite-type invariants of knots [V]. He completed the space
of all knots by adding degenerate ones and got an infinite-dimensional
contructible Euclidean space $E$. The degenerate knots form a hypersurface
$D$ (discriminant) in this space. Isotopy classes of knots are exactly
connected components of the complement to $D$. Alexander duality relates
$H^0(E \setminus  D)$ to the homology of $D$ which is then studied by using
its natural stratification.
   One would like to apply this approach to the classification of 3-manifolds.
This is the subject of the present paper.
  We construct our space $E$ by a version of Pontryagin-Thom construction.
More precisely $E$ is the space of maps $f: {\bf R}^8 \rightarrow { \bf  R}^5$,
such that
$f$, when  restricted  to the complement of the ball of sufficiently large
 but not fixed radius is just a projection
${\bf R}^8\rightarrow{ \bf R}^5$. This is obviously an affine space,
since $f_1+(1-t)f_2$ is also in $E$  if  $f_1$  and $f_2$ are.
The discriminant $D$ consists of maps $f: {\bf R}^8 \rightarrow{ \bf  R}^5$
which have
 critical point with critical value 0. If $f \notin D$, then $f^{-1}(0)$ is a
smooth
punctured 3-manifold, i.e., a compact manifold from which one point is deleted.
Our main result is as follows:

\proclaim Theorem 1. To each connected  component of $E- D$
there corresponds a   homeomorphism class of 3-dimensional spin manifolds.
For any connected spin  manifold there are exactly 4 connected components of
$E-D$ giving its homeomorphism type. In this case the
gauge group $\pi_{8}(SO(5))={ \bf Z_2}\oplus{ \bf Z_2}$, which  acts on $E-D$,
 permutes the four chambers corresponding to our manifold.

By a spin manifold we understand a pair $(M,\theta)$ where  $M$ is an oriented
3-manifold, and  $\theta$ is a spin structure on M. Two spin manifolds
$(M,\theta)$
and $(M',\theta')$ are called homeomorphic, if there exists a homeomorphism
 $M\rightarrow M'$ taking $\theta$ to $\theta'$.

The construction is as follows: given a connected punctured manifold $ \dot
M^3$ there
is unique isotopy type of embedding $ \dot  M^3 \hookrightarrow {\bf R}^8$.  A
representation of $M$ as $f^{-1}(0)$ as above gives a framing of $M$, i.e. a
trivialization of its normal bundle in ${\bf R}^8$. The set of such framings is
found to be ${ \bf  Z}\oplus\Sigma(M)$ for connected nonpunctured manifold M,
where
 $\Sigma(M)$ is the set of spin structures on $M$. By puncturing manifold
 we kill the  {\bf Z}-component of the framing and are left only with
$\Sigma(M)$, which depends on the surgery presentation of manifold $M$.
Given a framing, the question of representing $M$ as
$f^{-1}(0)$ with $df$ inducing this framing, is a certain problem of
obstruction theory.
 It turns out that this problem has a nice answer:
\proclaim Theorem 2. Let $M$ be a connected 3-manifold.
 For any spin structure $\theta$ on $M$  there is
 $f:R^8\rightarrow R^5$ as above such that
$f\notin D$ and further,  $M=f^{-1}(0)$ and $df$ induces the
 spin structure $\theta$.
 The set of isotopy classes of such $f$ consists of 4 elements.

\vskip .3cm

{ \bf Acknowledgements}. I want to thank my advisor S. Weinberger who
supervised this work,
O.Viro who attracted my attention to the importance of spin structures and also
S. Ackbulut,
R.Kirby, G.Kuperberg, M. Polyak for useful discussions.

 \vskip 1cm

\vfill\eject

\centerline { \bf  1.  Proof of Theorem 1.}

\vskip .7cm

   First we introduce some  definitions.
Let's denote by $E$ the infinite-dimensional affine space of 5-tuples
of functions $f=(f_1,....,f_5)$ on ${\bf R}^8$ such that outside a ball of
sufficient large radius we have $f_i(x_1, ..., x_8) = x_i$ for $i=1, ..., 5$.

\proclaim Definition 1. Denote by $D$ and call the  discriminant the subset in
$E$ which consists of such tuples $f$ that on the submanifold $\{f=0\}$ there
 exists a point, s.t. the Jacobian $Jf$ degenerates.

 \proclaim   Definition 2. Connected components of $E-D$ we will call chambers.

 \proclaim Definition 3.  Part of the discriminant separating two chambers is
 called a wall.

  For each chamber V and any $f=(f_1,.....f_5)$ in V the set of zeros
$M_f=\{f=0\}$ is a smooth submanifold in $R^8$.

  Since by Milnor's theorem any  connected 3-manifold can be obtained from the
other by
the surgery on the sequence of knots, we can assume that to each wall there
assigned a knot, over which we have to do surgery to pass to another chamber,
containing manifold of different homeomorphism type.
 For example there is a chamber which contains $S^3$. Walls of this chamber
correspond to different knots, which are already combinatorially classified
by Vassiliev.
 To pass from a chamber giving a manifold with
 two connected components to  one that
contains the connected sum of them we need not just (2,2)-type surgery as above
but
also (3,1)-type.
   In the case of 3-manifolds all three categories PL, TOP, Diff are
equivalent [Mo]. For convenience we will work with smooth manifolds. Note
that since all 3-manifolds are spin,  their tangent bundle is trivial.

\proclaim Lemma 1. The normal bundle to $M^3$ in $S^8$ is trivial.

Proof. Let $r$ denote the 1-dimensional trivial bundle.
 Note that the normal
bundle to $M^3$ in $S^8$ is stably trivial since characteristic
classes $w_2$, $\tau_1$ vanish. It is trivial because having a $k$-bundle $E$
overan  $n$-complex $(n< k+1)$ such that  $E\oplus r^l$ is trivial  we can
prove by induction that $E\oplus r$ is trivial.
Now,
$E\oplus r$ is isomorphic to $r^k+1$, so  we get a map from $M$ into
the projective space
$P^k$ defined by the position of the fiber of the
summand $r$. Since $n<k+1$,
 the image of a line can be homotoped to $Id$. Thus $E$ is trivial
for $n<k+1$, which holds in our case.

\vskip .3cm

Now we extend the (trivial) normal
${ \bf  R}^5$-bundle over ${\bf R}^8$ and consider
our  punctured 3-manifold $M$ as the locus of common zeroes
of five functions on $R^8$.
It is always possible to obtain $M^3$ in this way in some tubular
neighborhood $\nu M^3$ of $M^3$ in $R^8$ which we identify with the normal
bundle.
To do this, we should
just consider the tautological bundle whose fiber over a point $(x,v)$
of $\nu M^3$ is the fiber of $\nu M^3$ over $x$.
In $\nu M^3$ our manifold $M^3$ can be represented as the zero
section of this bundle.

We want to extend this section  (i.e., a map $f: \nu \dot  M^3\rightarrow
{ \bf  R}^5 - \{0\} \sim S^4$) to a map of all $R^8$ which does
not vanish outside $M^3$. By the standard obstruction theory
[D-N-F] the obstruction to  extension to the $i$th skeleton  in our setting
lies in
$$H^i(B(R)/ \nu \dot  M^3,d\nu( \dot  M^3) \cup (S^7/\nu(S^2)),
\pi_{i-1}(S^4))=$$
$$=H_{8-i}(B(R)/\nu \dot M^3; \pi_{i-1}(S^4)=H_{8-i}(S^8/M^3;\pi{i-1}(S^4)$$
 These groups can be easily calculated:
$$\cases{ H^{i-1}(M^3) = 0, i>4\cr
H^i(M^3, \pi_{i-1}(S^4)) = 0, i \leq 4.}$$
So the only possible obstruction lies in the group corresponding to i=0:
${ \cal O}\in { \bf Z}_{12}$.

 The number of extensions is given by the last
homology group:
$$H^8(S^8/M^3, \pi_8(S^4)) = { \bf  Z_2}\oplus { \bf  Z_2}.$$

\proclaim Lemma 2. We have an identification
$[M: SO(5)] \simeq { \bf  Z} \oplus H^1(M, { \bf  Z}_2)$.

\noindent {\sl Proof.} Consider the fibration
$${\rm Spin}(5) \rightarrow SO(5) \rightarrow B{ \bf  Z}_2.$$
 From this fibration we find an exact sequence
$$[M: {\rm Spin}(5)] \rightarrow  [M:SO(5)] \rightarrow
[M: B{ \bf  Z}_2],$$
so
$$0\rightarrow { \bf  Z} {\leftarrow \atop \rightarrow}
[M: SO(5)] \rightarrow H^1(M, { \bf  Z}_2) \rightarrow 0.$$
This sequence splits by the map $d$, which takes the
 map
$M\rightarrow SO(5)$ into its degree (defined in virtue of
equality $\pi_3(SO(5)) = { \bf  Z}$).

\vskip .2cm

We can write, therefore, our map obstruction for punctured manifold $ \dot
M^3$
in the form
$$\phi:  H^1(M, { \bf  Z}_2) \rightarrow
{ \bf  Z}_{12}.$$
We want to show that every element $\theta \in H^1(M, { \bf  Z}_2)$
is mapped by $\phi$ to 0.

\proclaim Lemma 3. For any 3-manifold $M$ and any spin structure
$\theta$ on $M$ obstruction $ \cal  O$ vanishes
for $(M, \theta)$.

\noindent {\sl Proof.}  From [Wa] we know that any 3-manifold embeds in
5-space and if  $(M,fr)$ represents zero in $\pi_3^{st}$ then there exists
embedding in $R^5$ realizing this framing (i.e. $\pi_3^{st}={ \bf Z}_{24}$
is the  obstruction group to getting framed 0-cobordism in $S^5$).
By Rohlin [R] we know that the third spin cobordism group is trivial, i.e.
each $M^3$ with a given spin structure is spin cobordant to $S^3$.
Now we have {\bf Z} possibilities to extend a given spin structure to a framing
on $M^3$ and those which represent zero in $\pi_3^{st}$ will provide a framed
0-cobordism. Thus each spin structure can be extended to a framing in
such a way that $(M, {\rm}fr)$ is given by equations in $S^5$.
This construction can be raised to $S^8$ in the following way.
  Suppose $M^3$ is  a complete intersection in $S^5$ and $W_1^4$ and $W_2^4$
are 4-manifolds whose intersection is $M^3$. In $S^5$ any closed 4-manifold
bounds. Let $V_1$ and $V_2$ be 5-manifolds whose boundaries are $W_1$
and $W_2$. Now consider $S^5$ as a submanifold of $S^6$ and construct
"pushouts" $V_1^+$, $V_1^-$, s.t. their common boundary is $W_1$. Make
the same construction for $V_2$. Construct $V_1^+\cup V_1^- = X_1^5$ and
$V_2^+ \cup V_2^- = X_2^5$.
 Now $M = S^5 \cap X_1^5 \cap X_2^5$. So the construction is
raised to $S^6$. Then raise it to $S^8$ by induction and to ${\bf R}^8$ by
puncturing
sphere.
The obstruction that we calculated measures how much our punctured spin
 manifold differs from being a complete intersection.
 Thus  obstruction  vanishes for all spin structures on M.

\proclaim Lemma 4. The action of the gauge group $\pi_8(SO(5))$  on the space
of 3-manifolds
permutes, in a simply transitive way, the  chambers corresponding to any given
connected spin manifold.

\noindent {\sl Proof.} The solution of the obstruction theory problem
provided us with four  different extentions for every  spin 3-manifold, which
are given
by the last homotopy group $\pi_{8}(S^4)={ \bf Z_2}\oplus{ \bf Z_2}$.
This has a natural geometrical explanation. Different trivializations of
5-bundle over ${\bf R}^8$  are classified by
 $\pi_{8}(SO(5))={ \bf Z_2}\oplus{
\bf Z_2}$.
 This is a gauge
group acting in 5-bundle and the isomorphism
 $$\pi_{8}(SO(5))\rightarrow\pi_{8}(S^4)$$
is induced by the map of the spaces. $SO(5)$ acts on $S^4=R^5-\{0\}$ and those
maps
induce the action of a gauge group on the set of all extensions, since all
other
obstruction groups are trivial.
Thus the quotient of the space of 3-manifolds by the action of the gauge group
gives unique extension,i.e., chamber for any spin 3-manifold.

\vfill\eject

\centerline { \bf  2. Axiomatics.}

\vskip 1cm

 It is apriori not obvious that one cannot
write  Vassiliev-Ohtsuki-type axioms for arbitrary 3-manifolds without any
additional structure.  However, the following argument which belongs to O.Viro
shows that
in this situation all invariants of finite type will be trivial.

Suppose we have an invariant $\alpha$ of order $n$, i.e., its $n+1$st
``Vassiliev derivative" is zero. Then  we show that it should be constant
(for $M^3$ without additional structure).

Indeed, take the last wall over which we are taking the derivative,
and write the corresponding difference.
$$W_{L_1,...,L_{n-1}}^M - W_{L_n}^{M'=M_{L_n}} = 0. $$
Thus the derivative over the last wall equals the derivative over all previous
walls
and doesn't depend on it. Since everything is invariant under Kirby
moves, adding a small unknotted component to the link in the last wall
won't change the invariant (but will change the link). Doing it inductively we
can change
the surgery coefficient on each component of the link and make the
components unlinked (by application of Kirby moves). By
a sequence of such moves we can get from any 3-manifold $M$ to
$S^3$, simplifying the link. This will imply that $\alpha(M)=
\alpha(S^3)$, so $\alpha$ is constant.
We show that for manifolds with spin structures the above argument
doesn't work.

In the case of manifolds with spin structures Kirby moves were described in
[K-M]. Suppose $M$ is obtained from $S^3$ by surgery on a link
$L$. Suppose $C$ is a characteristic sublink defining the
spin structure. Recall that $C$ is called characteristic,
if $(C,L_i) \cong (L_i, L_i) ({\rm mod} 2)$ for any $L_i\i L$.
The moves are as follows.

\item{(1)} Add or delete a disjoint unknotted component $K$ with
framing $\pm 1$, and set $C' = C+K$.

\item{(2)} If $i\neq j$, slide $L_i$ over $L_j$, i.e., set
$L_i' = L_i + L_j$ and
$$C' = \cases{C, \quad L_i\notin C\cr
C-(L_i+L_j) + L'_i, \quad L_i, L_j \notin C\cr
C-L_i + (L_j+L'_i), \quad L_i\i C, L_j \notin C.}$$

$$C'=\cases {C+K, \quad (C,K) \,\, {\rm even}\cr
C,\quad (C,K) \,\, {\rm odd}}$$

These moves are restrictive comparing to original Kirby moves. For example
if we have two characterstic strands  with $+$ -crossing and we put an
unknotted component over them to change sign of the crossing, then we cannot
take this component off the link after performing surgery since it becomes
characteristic.

Thus the previous argument (for manifolds without spin structures)
which implied triviality of all Vassiliev invariants, does not
work in the spin case.

\vskip .3cm

 In 1985 A. Casson introduced a new invariant
$\lambda$ of
oriented homology 3-spheres. It has two different descriptions.
One is as the intersection number of two Lagrangian subvarieties
parametrizing representations of the fundamental groups of handlebodies
in a Heegard decomposition of $M$, the intersection takes in the
symplectic variety of representations of the fundamental group
of the Riemann surface, the common boundary of the two handlebodies.

The other definition of the Casson invariant is by the surgery
formula:
$$ \lambda(S^3) = 0, \quad \lambda(M(K_n)) - \lambda(M(K_{n-1})) =
(1/2) \Delta''_K(1),$$
where $K$ is a knot in a homology 3-sphere, $M(K_n)$ is the
$(1,n)$th Dehn surgery on $K$, and $\Delta_K(t)$ is the Alexander
polynomial of $K$.

This formula suggests that Casson's invariant should be a Vassiliev
type invariant of homology spheres: the difference of values
of $\lambda$ for $M(K_n)$, $M(K_{n-1})$ (Vassiliev's discrete
derivative) is expressed in terms of the knot $K$, which
corresponds to the wall in the space of 3-manifolds separating
$M(K_n)$ and $M(K_{n-1})$. This suggests
definition which
was  given by S. Garoufalidis [G] and H. Ohtsuki [O].

In the case of homologically nontrivial 3-manifolds because of the argument
above
nontrivial Vassiliev theory can be built only for manifolds with additional
structure.
The construction of the space of 3-manifolds suggests the only natural
axiomatics in this situation:

\proclaim Definition.  A  map $v:\{{\rm homeo}\quad {\rm types} \quad
{\rm of} \quad {\rm  spin} \quad {\rm 3-manifolds}\} \rightarrow
C$ is called a finite type invariant of at most order k if it satisfies the
condition:
$$\sum_{{\rm char} L' \i L}(-1)^{\# L'}v(M_{L'})=0\leqno (1)$$
where L' is a characteristic sublink of L,
as well as
 the following
axioms:
$$I( \dot  M)=I(M). \leqno (2)$$
$$I(M_1\# M_2)=I(M_1)+I(M_2).\leqno (3)$$
$$I(S^3)=0.\leqno (4)$$

The above axiomatics suggests that one should consider spinor modifications of
known invariants. It was shown that Rozanskii-Witten invariant [R-W] and
universal Vassiliev invariant introduced by Thang Le [Le] restrict to
C.Lescop's [L]
 generalization of Casson's invariant. This invariant vanishes for manifolds
with first Betti
number greater that four. Perhaps spinor modifications of this invariant or
it's perturbations will be nontrivial for manifolds with higher Betti numbers.

Now we introduce the first simple example of Vassiliev invariant of finite
order.
Start with $M^3$ with spin structure. By Rohlin's theorem all 3-manifolds
are spin cobordant. Consider Euler characteristic of spin 0-cobordism minus 1
modulo 2. Denote it $I(M,spin)$.

\proclaim Lemma 5.  $I(M,spin)$  is an invariant of spin manifold $M^3$.
It equals Rohlin's invariant modulo 2.

Proof. We know that if $M^3$  bounding $W^4$ is connected then $W^4$ can be
assumed to have only one 0-handle and some 2-handles [K]. Thus $I(M,spin)$
gives us the rank of second cohomology modulo 2.
  Since  closed spin 4-manifold has a unimodular, symmetric
and even intersection form, it's rank is even and our invariant is zero, thus
it is well-defined
 Note that the signature of the intersection matrix modulo 2 equals its rank
modulo 2.
Rohlin's invariant is defined as the signature of the intersection matrix
modulo 16.
Thus the constructed invariant $I(M,spin)$ equals Rohlin's invariant modulo 2.

\proclaim Lemma 6. Invariant $I(M,spin)$ is  finite type Vassiliev of order 1.

Proof. By making surgery over a knot we add a handle to spin cobordism, i.e.
invariant increases by 1, i.e. its derivative is constant.Taking alternated
sum over four chambers adjacent to selfintersection of discriminant of
codimension two we get zero. Thus second order Vassiliev derivative
is zero and we get  an order one Vassiliev invariant.

 It is important to understand which singular manifolds form the discriminant
D.
\proclaim Lemma 7. Discriminant D consists of 3-manifolds with Morse
singularities. Codimension n selfintersection of discriminant corresponds
to a manifold with n singular poins. Singular points have a type of  a cone
over torus for (2,2)-type surgeries and a cone over sphere for (3,1) and
(1,3)-type surgeries.

Proof. The whole picture can be described in terms of spin  4-cobordisms
connecting our 3-manifolds. Surgery over a knot corresponds to a  Morse
decomposition
of signature (2,2). It is well known how quadratic form of signature (2,2)
looks like. It has 2 generators and after projection to $RP^3$ is given by
quadratic
 homogeneous equation. So, it is a cone over the product of projective spaces,
i.e. a cone over a  torus.
 There are also two other surgeries of signature(1,3) and (3,1) (attaching
and eliminating handles). Those correspond to a cone over sphere.
 Passing from one chamber to another we are making Dehn surgery over a
knot which is assigned to the wall .

 Another interesting question which arises in connection with the space of
3 manifolds is  understanding of the topology of the chambers:

\proclaim Conjecture. The homotopy type of the chamber corresponding to $M$ is
\hfill\break
$K(\pi_0(Diff(M)),1)$.

Several cases are already known. For example for $S^3$ the diffeomorphism group
is homotopy equivalent to  $O(4)$. Also,
$Diff(S^1\times S^2)=O(2)\times O(3)\times\Omega SO(3)$, see [H1].
 For M Haken
$\pi_0 Diff(M)=Out(\pi_1(M)$ and the components of Diff(M) are usually
contractible. The exception is when M is Seifert-fibred, with all fibers
coherently orientable, and in this case the components of Diff(M) are
homotopy equivalent to $S^1$, [H2].
  Analogous theorem in the case of knots have been recently proved by A.
Hatcher.

\vskip 1cm

\vfill\eject
\centerline { \bf  References}
\vskip 1cm

[AM] S. Akbulut, J.C. McCarthy. Casson's invariant for oriented homology
3-spheres: an exposition. Princeton Math Notes, Princeton, 1990.

[B-N] D. Bar-Natan, On the Vassiliev knot invariants,

[B-L] J. Birman, X-S Lin, Knot polynomials and Vassiliev invariants, Invent.
Math. 111,1993,pp.225-270.

[C] C. Curtis,Generalized Casson invariants for SO(3), U(2), Spin(4) and
SO(4), preprint to appear in Trans AMS, 1993.

[D-N-F] B. Dubrovin, S. Novikov, A. Fomenko, Modern Geometry, Nauka, 1984.

[F] A. Floer, An instanton invariant for 3-manifolds, Commun. Math. Phys. 118,
1988,pp.215-240.

[G] S. Garoufalidis, On finite type 3-manifold invariants, preprint 1994.

[G-M] L. Gillou, A. Marin, A la recherche de la Topologie Perdu, Birkhauser,
1986.

[H1] A. Hatcher,Proc.A.M.S. 83,1981, pp.427-430.

[H2] A. Hatcher, Topology 15, 1987, pp. 343-348.

[L] C. Lescop, On the Casson-Walker invariant: global surgery formula and
extension to closed oriented 3-manifolds, preprint 1993.

[K] R. Kirby, The topology of 4-manifolds, Springer-Verlag

[K-M] R. Kirby, P. Melvin, The 3-manifold invariants of Reshetikhin-Turaev
for $sl(2,{ \bf  C})$, Inventiones, 105,1991,pp. 473-545.

[Ko] M. Kontsevich, Vassiliev's knot invariants, Adv. in Sov. Math., 16(2),
1993, pp.137-150.

[Le] Thang Le, Universal Vassiliev invariant, preprint 1996.

[Mu] H. Murakami, Quantum $SO(3)$ invariants dominate Casson's $SU(2)$
invariant , preprint 1993.

[O] T. Ohtsuki, Finite type invariants of integral homology 3-spheres,
preprint 1993.

[R] V. Rohlin, Classification of the mappings of (n+3)-dimensional sphere,
Dokl. Akad. Nauk, 1951, \# 4, pp. 541-544.

[Ro] D. Rolfsen, Knots and links,Publish or Perish, 1976.

[R,W] L. Rozansky, E. Witten, Hyper - Kahler Geometry and Invariants
of 3-manifolds, preprint 1996.

[V] V. A. Vassiliev, Complements of discriminants of smooth maps, Trans. of
Math. Mono. 98, Amer. Math. Soc., Providence, 1992.

[W] K. Walker, An extension of Casson's invariant, Princeton Univ. Press,
1993.

[Wa] C.T.C. Wall, All 3-manifolds imbed in 5-space.

\bye